\documentclass[11pt]{amsart}
\usepackage[utf8]{inputenc}
\usepackage[T1]{fontenc}
\usepackage{lmodern}
\usepackage{amsmath,amssymb,amsthm,mathtools}
\usepackage{booktabs,array}
\usepackage{tikz-cd}
\usepackage[margin=2.6cm]{geometry}
\usepackage{microtype}
\usepackage[hidelinks]{hyperref}

\newtheorem{theorem}{Theorem}[section]
\newtheorem{lemma}[theorem]{Lemma}
\newtheorem{proposition}[theorem]{Proposition}
\newtheorem{corollary}[theorem]{Corollary}
\theoremstyle{definition}
\newtheorem{definition}[theorem]{Definition}
\theoremstyle{remark}
\newtheorem{remark}[theorem]{Remark}

\newcommand{\C}{\mathbb{C}}
\newcommand{\Pj}{\mathbb{P}}
\newcommand{\lt}{\dashv}
\newcommand{\rt}{\vdash}
\newcommand{\Ann}{\operatorname{Ann}}
\newcommand{\Aut}{\operatorname{Aut}}
\newcommand{\Der}{\operatorname{Der}}
\newcommand{\GL}{\mathrm{GL}}
\newcommand{\SO}{\mathrm{SO}}
\newcommand{\Ot}{\mathrm{O}}
\newcommand{\Dl}[1]{\Delta^{\lt}_{#1}}
\newcommand{\Dr}[1]{\Delta^{\rt}_{#1}}
\newcommand{\N}[2]{\mathcal{N}_2[#1:#2]}
\newcommand{\Nz}{\mathcal{N}_2^0}
\newcommand{\diag}{\operatorname{diag}}

\begin{document}

\title[Three-dimensional nilpotent LSDA]{The algebraic classification of three-dimensional nilpotent left-symmetric dialgebras}

\author{Ikrom M. Rikhsiboev}


\maketitle

\begin{center}\small
Westminster International University in Tashkent, Tashkent, Uzbekistan\\
E-mail: \href{mailto:ikromr@gmail.com}{ikromr@gmail.com}
\end{center}

\begin{abstract}
Left-symmetric dialgebras (LSDA) are to left-symmetric algebras what Leibniz algebras are to Lie algebras. We classify the three-dimensional nilpotent LSDA over the complex numbers. To do this we adapt the annihilator-extension (Skjelbred--Sund) method to dialgebras, using the annihilator with respect to both products on both sides, and we compute the second cohomology of all two-dimensional nilpotent LSDA. Up to isomorphism, besides the null and split algebras there are nine non-split algebras and families; the $2$-step nilpotent families depend on four and three parameters. We give isomorphism criteria and separating invariants, identify the left-symmetric and diassociative algebras in the list, and complete the two-dimensional classification of an earlier paper.
\end{abstract}

\section{Introduction}

An algebra $(A,\cdot)$ is left-symmetric if the associator $(x,y,z)=(xy)z-x(yz)$ is symmetric in $x$ and $y$. These algebras go back to Cayley's work on rooted trees \cite{Cayley} and arose independently in Vinberg's theory of homogeneous convex cones \cite{Vinberg} and in the study of affine structures on Lie groups. They are also called pre-Lie, Vinberg or Koszul--Vinberg algebras. The commutator of a left-symmetric product is a Lie bracket, and left-symmetric algebras appear in operad theory \cite{ChapotonLivernet}, deformation theory, renormalisation and Rota--Baxter theory \cite{LiHouBai}. Burde's survey \cite{Burde} collects their main applications in geometry and physics. Classifications are known in low dimensions: Kim \cite{Kim} described the complete left-symmetric structures on nilpotent Lie algebras of small dimension, and Bai \cite{Bai} classified the three-dimensional left-symmetric algebras over $\C$.

Loday introduced dialgebras, vector spaces with two products $\lt$ and $\rt$, as the natural setting for Leibniz algebras \cite{Loday93,LodayEtAl}: an associative dialgebra (diassociative algebra) gives a Leibniz algebra in the same way an associative algebra gives a Lie algebra. Diassociative algebras have since been studied structurally \cite{Gonzalez,LinZhang}, and the nilpotent ones were classified in dimension three \cite{RRB10} and four \cite{RRB15}. Felipe \cite{Felipe} defined \emph{left-symmetric dialgebras} (LSDA) as the dialgebra analogue of left-symmetric algebras: they are to left-symmetric algebras what Leibniz algebras are to Lie algebras. When $\lt=\rt$ the axioms reduce to the left-symmetric identity, and every diassociative algebra is an LSDA. The author began the structure theory of LSDA in \cite{R16}, introducing annihilators and nilpotency and classifying the two-dimensional LSDA over a field of characteristic zero.

The standard tool for nilpotent algebras of small dimension is the method of annihilator (central) extensions of Skjelbred and Sund \cite{SkjelbredSund}, first used for nilpotent Lie algebras \cite{deGraaf,GozeKhakimdjanov}. Every nilpotent algebra is an annihilator extension of a smaller one, and isomorphism classes of extensions correspond to orbits of the automorphism group on subspaces of a second cohomology space. Hegazi, Abdelwahab and Calder\'on Mart\'in adapted it to Malcev and other non-associative algebras \cite{HAC}, and it has since been applied to many varieties (Leibniz \cite{AOR}, Jordan, Novikov, bicommutative, Zinbiel and others); see the survey by Kaygorodov \cite{Kaygorodov}. Boza, Fedriani, N\'u\~nez and Tenorio \cite{BFNT} review the history of Lie algebra classifications.

For LSDA no classification beyond dimension two is known, and the two-dimensional list of \cite[Theorem~2.6]{R16} needs to be completed (Remark~\ref{rem:2016}). The cohomology needed for annihilator extensions has not been set up for LSDA either. A difficulty specific to dialgebras is the choice of annihilator: the extension method needs the subspace killed on both sides by both products, which is smaller than the ideal used in \cite[Definition~2.3]{R16}.

In this paper we classify all three-dimensional nilpotent left-symmetric dialgebras over $\C$ (Theorem~\ref{thm:main}), using annihilator extensions adapted to dialgebras. Besides the null and split algebras, there are nine non-split algebras and families, two of which depend on four and three parameters.

\section{Preliminaries}

All vector spaces are finite-dimensional over $\C$. For subspaces $U,W$ of an algebra and a product $\circ$, $U\circ W$ is the span of all $u\circ w$.

\begin{definition}[\cite{Felipe,R16}]
A \emph{left-symmetric dialgebra} (LSDA) is a vector space $S$ with two bilinear products $\lt$ and $\rt$ such that, for all $x,y,z\in S$,
\begin{align}
x\lt(y\lt z)&=x\lt(y\rt z),\tag{1.1}\\
(x\rt y)\rt z&=(x\lt y)\rt z,\tag{1.2}\\
x\lt(y\lt z)-(x\lt y)\lt z&=y\rt(x\lt z)-(y\rt x)\lt z,\tag{1.3}\\
x\rt(y\rt z)-(x\rt y)\rt z&=y\rt(x\rt z)-(y\rt x)\rt z.\tag{1.4}
\end{align}
We keep the numbering of \cite{R16}. A homomorphism of LSDA is a linear map preserving both products.
\end{definition}

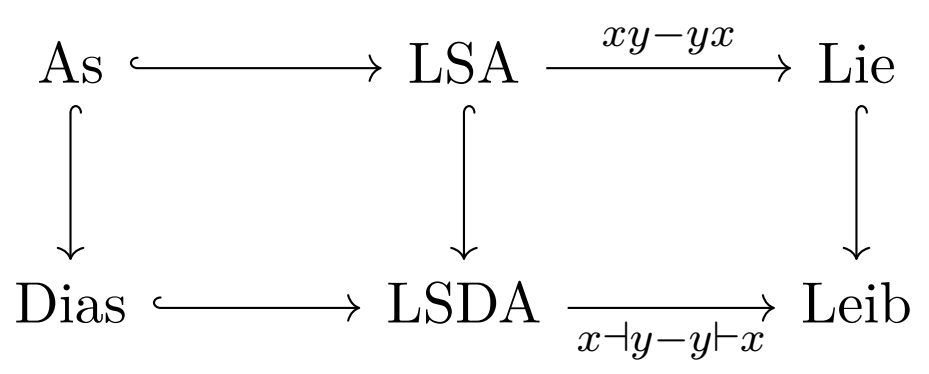
\begin{figure}[ht]
\[
\begin{tikzcd}[column sep=large,row sep=large]
\mathrm{As} \arrow[r, hook] \arrow[d, hook] & \mathrm{LSA} \arrow[r, "{xy-yx}"] \arrow[d, hook] & \mathrm{Lie} \arrow[d, hook]\\
\mathrm{Dias} \arrow[r, hook] & \mathrm{LSDA} \arrow[r, "{x\lt y-y\rt x}"'] & \mathrm{Leib}
\end{tikzcd}
\]
\caption{The extended Loday diagram: associative algebras (As), diassociative algebras (Dias), left-symmetric algebras (LSA), left-symmetric dialgebras (LSDA), Lie and Leibniz (Leib) algebras.}\label{fig:loday}
\end{figure}

\begin{remark}\label{rem:special}
If $\lt=\rt$, then (1.1) and (1.2) are trivial and (1.3), (1.4) both become the left-symmetric identity, so an LSDA with equal products is a left-symmetric algebra. Every diassociative algebra \cite{LodayEtAl} is an LSDA: its axioms make both sides of (1.3) and of (1.4) zero and give (1.1), (1.2) directly. Moreover, the commutator $[x,y]=x\lt y-y\rt x$ of an LSDA satisfies the (right) Leibniz identity $[[x,y],z]=[[x,z],y]+[x,[y,z]]$ \cite{Felipe}, \cite[Theorem~1.13]{R16}; for diassociative algebras this is Loday's functor $\mathrm{Dias}\to\mathrm{Leib}$, and for $\lt=\rt$ it is the commutator $\mathrm{LSA}\to\mathrm{Lie}$. These relations form the extended Loday diagram of \cite{R16} (Figure~\ref{fig:loday}), in which hooked arrows are inclusions of subcategories and the horizontal arrows on the right are the commutator functors.
\end{remark}

\begin{definition}\label{def:ann}
Let $S$ be an LSDA.
\begin{enumerate}
\item The \emph{annihilator} of $S$ is $\Ann(S)=\{x\in S: x\lt S=S\lt x=x\rt S=S\rt x=0\}$. It is an ideal of $S$.
\item Put $S^1=S$ and $S^k=\sum_{i+j=k}\bigl(S^i\lt S^j+S^i\rt S^j\bigr)$ for $k\ge 2$. Then $S^{k+1}\subseteq S^k$ for all $k\ge1$: this is clear for $k=1$, and for $k\ge2$ every summand $S^i\circ S^j$ ($\circ\in\{\lt,\rt\}$) of $S^{k+1}$ has $i\ge2$ or $j\ge2$, so by induction it lies in $S^{i-1}\circ S^j$ or $S^i\circ S^{j-1}$, which are summands of $S^k$. The algebra $S$ is \emph{nilpotent} if $S^t=0$ for some $t$; the least such $t$ is the \emph{nilpotency index}.
\item $S$ is \emph{split} if $S=S'\oplus\C e$, a direct sum of ideals with $e\in\Ann(S)$, and \emph{non-split} otherwise.
\end{enumerate}
\end{definition}

The annihilator in \cite[Definition~2.3]{R16} is the intersection of the right $\lt$-annihilator $\{x: S\lt x=0\}$ and the left $\rt$-annihilator $\{x: x\rt S=0\}$. It contains $\Ann(S)$ and may be strictly larger. Annihilator extensions need the subspace killed by both products on both sides, which is $\Ann(S)$ above.

\begin{lemma}\label{lem:ann}
Every non-zero nilpotent LSDA $S$ has $\Ann(S)\neq 0$. More precisely, if $t$ is the nilpotency index of $S$, then $0\neq S^{t-1}\subseteq\Ann(S)$.
\end{lemma}

\begin{proof}
Since $S\neq0$ we have $t\ge 2$, and $S^{t-1}\neq 0$ by minimality of $t$. The terms with $(i,j)=(1,t-1)$ and $(t-1,1)$ in the definition of $S^t$ give $S\lt S^{t-1}$, $S\rt S^{t-1}$, $S^{t-1}\lt S$, $S^{t-1}\rt S\subseteq S^t=0$. Hence $S^{t-1}\subseteq\Ann(S)$.
\end{proof}

Since $\Ann(S)$ is an ideal, $S/\Ann(S)$ is again a nilpotent LSDA, of smaller dimension by Lemma~\ref{lem:ann}. This is the basis of the inductive method below.

\begin{definition}\label{def:coh}
Let $A$ be an LSDA and $V$ a vector space. A pair $\theta=(\theta_\lt,\theta_\rt)$ of bilinear maps $A\times A\to V$ is a \emph{$2$-cocycle} if, for all $x,y,z\in A$,
\begin{align*}
&\theta_\lt(x,y\lt z)=\theta_\lt(x,y\rt z),\qquad \theta_\rt(x\rt y,z)=\theta_\rt(x\lt y,z),\\
&\theta_\lt(x,y\lt z)-\theta_\lt(x\lt y,z)=\theta_\rt(y,x\lt z)-\theta_\lt(y\rt x,z),\\
&\theta_\rt(x,y\rt z)-\theta_\rt(x\rt y,z)=\theta_\rt(y,x\rt z)-\theta_\rt(y\rt x,z).
\end{align*}
The cocycles form a space $Z^2(A,V)$. For a linear map $f\colon A\to V$, the pair $\delta f=\bigl(f(x\lt y),f(x\rt y)\bigr)$ is a \emph{coboundary}; the coboundaries form $B^2(A,V)$, and $H^2(A,V)=Z^2(A,V)/B^2(A,V)$. We write $\Ann(\theta)=\{x\in A:\theta_\circ(x,A)=\theta_\circ(A,x)=0 \text{ for } \circ\in\{\lt,\rt\}\}$. For $\theta\in Z^2(A,V)$ the \emph{annihilator extension} $A_\theta=A\oplus V$ has the products
\[
(x+u)\lt(y+w)=x\lt_A y+\theta_\lt(x,y),\qquad (x+u)\rt(y+w)=x\rt_A y+\theta_\rt(x,y).
\]
\end{definition}

\begin{lemma}\label{lem:cocycle}
Let $\theta$ be a pair of bilinear maps $A\times A\to V$. Then $A_\theta$ is an LSDA if and only if $\theta\in Z^2(A,V)$. Moreover $B^2(A,V)\subseteq Z^2(A,V)$.
\end{lemma}

\begin{proof}
Products involving $V$ vanish, so it suffices to check (1.1)--(1.4) on $x,y,z\in A$. For example $x\lt(y\lt z)=x\lt_A(y\lt_A z)+\theta_\lt(x,y\lt_A z)$ in $A_\theta$. Each identity in $A_\theta$ therefore splits into its $A$-component, which is the same identity in $A$, and its $V$-component, which is the corresponding condition of Definition~\ref{def:coh}. For $\delta f$, each condition is $f$ applied to an identity of $A$, so it holds.
\end{proof}

\begin{lemma}\label{lem:annext}
$\Ann(A_\theta)=\bigl(\Ann(A)\cap\Ann(\theta)\bigr)\oplus V$.
\end{lemma}

\begin{proof}
$V\subseteq\Ann(A_\theta)$. For $x,y\in A$ we have $x\lt y=x\lt_A y+\theta_\lt(x,y)$, with the two terms in $A$ and $V$. So $x\lt A_\theta=0$ iff $x\lt_A A=0$ and $\theta_\lt(x,A)=0$; the same holds for the other three conditions.
\end{proof}

\begin{lemma}\label{lem:every}
Let $S$ be an LSDA with $V=\Ann(S)\neq0$ and $A=S/V$. Then $S\cong A_\theta$ for some $\theta\in Z^2(A,V)$ with $\Ann(A)\cap\Ann(\theta)=0$.
\end{lemma}

\begin{proof}
Choose a complement $S=A'\oplus V$ and let $\pi$ be the projection onto $A'$ along $V$. Identify $A'$ with $A=S/V$; since $V$ is an ideal, the products of $A$ are $x\circ_A y=\pi(x\circ y)$. Put $\theta_\circ(x,y)=x\circ y-\pi(x\circ y)\in V$ for $\circ\in\{\lt,\rt\}$. Then $S=A_\theta$, so $\theta$ is a cocycle by Lemma~\ref{lem:cocycle}, and Lemma~\ref{lem:annext} gives $V=\Ann(S)=(\Ann(A)\cap\Ann(\theta))\oplus V$.
\end{proof}

Fix a basis $e_1,\dots,e_s$ of $V$ and write $\theta=\sum\theta_ie_i$ with $\theta_i\in Z^2(A,\C)$. For $\varphi\in\Aut(A)$ put $(\varphi\theta)(x,y)=\theta(\varphi x,\varphi y)$. Since $\Ann(A)\cap\Ann(\theta+\delta f)=\Ann(A)\cap\Ann(\theta)$, the condition $\Ann(A)\cap\Ann(\theta)=0$ depends only on the subspace $\langle[\theta_1],\dots,[\theta_s]\rangle\subseteq H^2(A,\C)$.

\begin{lemma}\label{lem:iso}
Let $\theta,\eta\in Z^2(A,V)$ with $\Ann(A)\cap\Ann(\theta)=\Ann(A)\cap\Ann(\eta)=0$.
\begin{enumerate}
\item $A_\theta$ is split if and only if $[\theta_1],\dots,[\theta_s]$ are linearly dependent in $H^2(A,\C)$.
\item If $A_\theta$ and $A_\eta$ are non-split, then $A_\theta\cong A_\eta$ if and only if there is $\varphi\in\Aut(A)$ with $\langle[\varphi\theta_1],\dots,[\varphi\theta_s]\rangle=\langle[\eta_1],\dots,[\eta_s]\rangle$.
\end{enumerate}
\end{lemma}

\begin{proof}
(1) If the classes are dependent, change the basis of $V$ so that $\theta_s=\delta f$ with $f\colon A\to\C$. Replacing $A$ by its graph $\{x-f(x)e_s\}$ changes $\theta$ by a coboundary and makes $\theta_s=0$; then $\C e_s$ is an ideal with complement $A\oplus\langle e_1,\dots,e_{s-1}\rangle$, so $A_\theta$ is split. Conversely, let $A_\theta=S'\oplus\C e$ be a direct sum of ideals with $e\in\Ann(A_\theta)=V$. Then $V=(V\cap S')\oplus\C e$. Choose a complement $S'=A''\oplus(V\cap S')$; all products lie in $S'$, so with respect to $A_\theta=A''\oplus V$ and a basis of $V$ ending with $e$, the last component of the cocycle vanishes. A change of complement and of basis of $V$ replaces the $\theta_i$ by invertible combinations plus coboundaries, so the classes $[\theta_i]$ are dependent.

(2) Let $\psi\colon A_\theta\to A_\eta$ be an isomorphism. By Lemma~\ref{lem:annext} both annihilators equal $V$, so $\psi(V)=V$. Hence $\psi(x+v)=\bar\varphi(x)+g(x)+Mv$ with $\bar\varphi\in\Aut(A)$, $g\colon A\to V$ linear and $M\in\GL(V)$. Comparing the $V$-components of $\psi(x\circ y)$ and $\psi(x)\circ\psi(y)$ for $\circ\in\{\lt,\rt\}$ gives
\[
M\theta+\delta g=\bar\varphi\eta.
\]
Conversely, if $\bar\varphi$, $M$ and $g$ satisfy this equation, the same computation shows that $\psi(x+v)=\bar\varphi(x)+g(x)+Mv$ is an isomorphism. Applying $\bar\varphi^{-1}$ gives $M(\varphi\theta)+\delta(g\bar\varphi^{-1})=\eta$ with $\varphi=\bar\varphi^{-1}$. Thus $A_\theta\cong A_\eta$ iff $[\eta_1],\dots,[\eta_s]$ are the images of $[\varphi\theta_1],\dots,[\varphi\theta_s]$ under some invertible matrix $M$. For two independent families this is exactly equality of their spans.
\end{proof}

\begin{corollary}\label{cor:proc}
The non-split $n$-dimensional nilpotent LSDA with $s$-dimensional annihilator correspond bijectively to the pairs $(A,W)$, where $A$ runs over the $(n-s)$-dimensional nilpotent LSDA up to isomorphism and $W$ over the $\Aut(A)$-orbits of $s$-dimensional subspaces $W\subseteq H^2(A,\C)$ with $\Ann(A)\cap\Ann(W)=0$. Every split algebra is $B\oplus\C$ with $B$ nilpotent of dimension $n-1$.
\end{corollary}

\begin{proof}
Lemmas~\ref{lem:ann}, \ref{lem:every} and \ref{lem:iso}.
\end{proof}

\section{Two-dimensional nilpotent LSDA and their cohomology}

\textbf{Notation.} In all multiplication tables only the non-zero products of basis vectors are listed. For an algebra with basis $e_1,\dots,e_n$ we write $\Dl{ij}$ (resp.\ $\Dr{ij}$) for the cocycle with $\theta_\lt(e_i,e_j)=1$ (resp.\ $\theta_\rt(e_i,e_j)=1$) and all other values $0$. For $\theta\in Z^2(A,\C)$ we put $\lambda_{ij}=\theta_\lt(e_i,e_j)$ and $\rho_{ij}=\theta_\rt(e_i,e_j)$. $\mathcal N_1$ and $\Nz$ denote the null algebras of dimensions $1$ and $2$.

\begin{proposition}\label{prop:dim2}
Every two-dimensional nilpotent LSDA over $\C$ is isomorphic to the null algebra $\Nz$ or to exactly one algebra
\[
\N{\alpha}{\beta}\colon\quad e_1\lt e_1=\alpha e_2,\quad e_1\rt e_1=\beta e_2,\qquad [\alpha:\beta]\in\Pj^1.
\]
Always $\Ann(\N{\alpha}{\beta})=\langle e_2\rangle$, and $\Aut(\N{\alpha}{\beta})$ consists of the maps $\varphi\colon e_1\mapsto pe_1+qe_2$, $e_2\mapsto p^2e_2$ with $p\neq0$.
\end{proposition}

\begin{proof}
Let $S$ be non-null. By Lemma~\ref{lem:ann}, $\dim\Ann(S)=1$, and $A=S/\Ann(S)$ is one-dimensional nilpotent, hence $A=\mathcal N_1$. Since $\mathcal N_1$ has zero products, every pair $\theta=(\alpha\Dl{11},\beta\Dr{11})$ is a cocycle and $B^2(\mathcal N_1)=0$, so $H^2(\mathcal N_1,\C)=\C^2$. The condition $\Ann(\mathcal N_1)\cap\Ann(\theta)=0$ means $\theta\neq0$. $\Aut(\mathcal N_1)=\C^*$ acts by $p^2$ on both coordinates, so the one-dimensional subspaces $\langle\theta\rangle$, that is the points $[\alpha:\beta]\in\Pj^1$, are pairwise inequivalent. By Corollary~\ref{cor:proc} this gives the list. An automorphism preserves $S^2=\langle e_2\rangle$, so $\varphi(e_2)=re_2$ and $\varphi(e_1)=pe_1+qe_2$; comparing $\varphi(e_1)\circ\varphi(e_1)$ with $\varphi(e_1\circ e_1)$ gives $r=p^2$.
\end{proof}

\begin{remark}[Corrections to \cite{R16}]\label{rem:2016}
In Case~2 of the proof of \cite[Theorem~2.6]{R16} ($\alpha_3=0$) the identities give $\alpha_1=0$ and $\alpha_2\beta_1=0$, but only the subcase $\alpha_2\neq0$, $\beta_1=0$, $\beta_2\neq0$ was treated. It yields the nilpotent algebra $e_2\lt e_2=ae_1$, $e_2\rt e_2=e_1$ ($a\neq0$), which is $\N{a}{1}$ after interchanging $e_1$ and $e_2$; this algebra does not appear in the list $L_2^1$--$L_2^6$ of the statement. The remaining subcases give three further algebras: $e_2\lt e_2=e_1$ (that is, $\N{1}{0}$), $e_2\rt e_2=e_1$ (that is, $\N{0}{1}$), and the non-nilpotent algebra $e_2\rt e_1=e_1$, obtained after the change $e_2\mapsto e_2-be_1$ removes a product $e_2\rt e_2=be_1$. The last algebra is not isomorphic to any of $L_2^1$--$L_2^6$, because its product $\lt$ is zero, while each of $L_2^1$--$L_2^6$ has $\lt\neq0$. Hence \cite[Theorem~2.6]{R16} should be supplemented by the nilpotent algebras $\N{\alpha}{\beta}$ of Proposition~\ref{prop:dim2} and by the algebra $e_2\rt e_1=e_1$; all of them satisfy (1.1)--(1.4).
\end{remark}

In what follows we normalise $\alpha=1$ when $\alpha\neq0$, so the algebras are $\N{1}{\beta}$ ($\beta\in\C$) and $\N{0}{1}$.

\begin{proposition}\label{prop:H2}
Let $A=\N{\alpha}{\beta}$. A pair $\theta$ is a cocycle if and only if
\[
\lambda_{22}=\rho_{22}=0,\qquad (\alpha-\beta)\lambda_{12}=(\alpha-\beta)\rho_{21}=0,\qquad \alpha\rho_{12}=\alpha\lambda_{12}+(\beta-\alpha)\lambda_{21},
\]
and $B^2(A,\C)=\langle\alpha\Dl{11}+\beta\Dr{11}\rangle$. Consequently:

\smallskip
\begin{center}\footnotesize
\begin{tabular}{@{}>{\raggedright}p{2.1cm}>{\raggedright}p{4.4cm}>{\raggedright}p{2.3cm}c>{\raggedright\arraybackslash}p{3.6cm}@{}}
\toprule
$A$ & Conditions on $\theta$ & $B^2$ & $\dim H^2$ & Basis of $H^2$\\
\midrule
$\mathcal N_1$ & none & $0$ & $2$ & $\Dl{11},\ \Dr{11}$\\
$\Nz$ & none & $0$ & $8$ & all $\Dl{ij},\ \Dr{ij}$\\
$\N{1}{\beta}$, $\beta\neq1$ & $\lambda_{12}=\lambda_{22}=\rho_{21}=\rho_{22}=0$, $\rho_{12}=(\beta-1)\lambda_{21}$ & $\langle\Dl{11}+\beta\Dr{11}\rangle$ & $2$ & $[\Dl{11}],\ [\Dl{21}+(\beta-1)\Dr{12}]$\\
$\N{0}{1}$ & $\lambda_{12}=\lambda_{21}=\lambda_{22}=\rho_{21}=\rho_{22}=0$ & $\langle\Dr{11}\rangle$ & $2$ & $[\Dl{11}],\ [\Dr{12}]$\\
$\N{1}{1}$ & $\lambda_{22}=\rho_{22}=0$, $\lambda_{12}=\rho_{12}$ & $\langle\Dl{11}+\Dr{11}\rangle$ & $4$ & $[\Dl{11}],[\Dl{21}],[\Dr{21}],[\Dl{12}+\Dr{12}]$\\
\bottomrule
\end{tabular}
\end{center}
\smallskip
Here $\dim Z^2$ equals $2,8,3,3,5$ in the five rows.
\end{proposition}

\begin{proof}
The only non-zero products of $A$ are $e_1\lt e_1=\alpha e_2$ and $e_1\rt e_1=\beta e_2$, so a term $\theta(u,v\circ w)$ or $\theta(u\circ v,w)$ of Definition~\ref{def:coh} is non-zero only when the inner pair is $(e_1,e_1)$. Evaluating the four conditions on all basis triples gives:
\begin{itemize}
\item first condition, $(e_i,e_1,e_1)$: $(\alpha-\beta)\lambda_{i2}=0$ for $i=1,2$;
\item second condition, $(e_1,e_1,e_j)$: $(\alpha-\beta)\rho_{2j}=0$ for $j=1,2$;
\item third condition, triples $(e_1,e_1,e_1)$, $(e_1,e_1,e_2)$, $(e_2,e_1,e_1)$, $(e_1,e_2,e_1)$: $\alpha(\lambda_{12}-\lambda_{21})=\alpha\rho_{12}-\beta\lambda_{21}$, $(\alpha-\beta)\lambda_{22}=0$, $\alpha\lambda_{22}=0$ and $\alpha\rho_{22}=0$;
\item fourth condition, triples $(e_2,e_1,e_1)$, $(e_1,e_2,e_1)$: $\beta\rho_{22}=0$.
\end{itemize}
All remaining triples give $0=0$. Since $(\alpha,\beta)\neq(0,0)$, these equations are equivalent to the stated ones. For $f=\mu e_1^*+\nu e_2^*$ we get $\delta f=\nu(\alpha\Dl{11}+\beta\Dr{11})$. Specialising to $\alpha=1,\beta\neq1$, to $[0:1]$ and to $[1:1]$ gives the rows of the table; for $\mathcal N_1$ and $\Nz$ every pair is a cocycle and $B^2=0$.
\end{proof}

\begin{lemma}\label{lem:action}
Let $A=\N{\alpha}{\beta}$, $\varphi\in\Aut(A)$ as in Proposition~\ref{prop:dim2}, and $k\in\C^*$ a rescaling of $V=\C$. Then $k\varphi\theta$ has
\begin{align*}
\lambda'_{12}&=kp^3\lambda_{12},\quad \lambda'_{21}=kp^3\lambda_{21},\quad \rho'_{12}=kp^3\rho_{12},\quad \rho'_{21}=kp^3\rho_{21},\\
\lambda'_{11}&=k\bigl(p^2\lambda_{11}+pq(\lambda_{12}+\lambda_{21})\bigr),\qquad \rho'_{11}=k\bigl(p^2\rho_{11}+pq(\rho_{12}+\rho_{21})\bigr).
\end{align*}
\end{lemma}

\begin{proof}
Expand $\theta(\varphi e_i,\varphi e_j)$ with $\varphi e_1=pe_1+qe_2$, $\varphi e_2=p^2e_2$ and use $\lambda_{22}=\rho_{22}=0$.
\end{proof}

\section{The main theorem}

\begin{theorem}[Main Theorem]\label{thm:main}
Every three-dimensional nilpotent left-symmetric dialgebra over $\C$ is isomorphic to exactly one of the following algebras: the null algebra $\mathcal N_3$, a split algebra $\N{\alpha}{\beta}\oplus\C$ with $[\alpha:\beta]\in\Pj^1$, or one of the nine non-split algebras and families $F_1$--$F_7$, $D$ and $T$ below. Only non-zero products are listed; the basis is chosen so that $e_3$ spans $\Ann(S)$ (for $F_7$, $e_2$ and $e_3$ span it), hence every product with $e_3$ as a factor is zero.

\smallskip\noindent\emph{Algebras with $S/\Ann(S)$ two-dimensional and non-null (nilpotency index $4$):}
\begin{itemize}
\item $F_1(\beta)$, $\beta\in\C\setminus\{1\}$: $e_1\lt e_1=e_2$, $e_1\rt e_1=\beta e_2$, $e_2\lt e_1=e_3$, $e_1\rt e_2=(\beta-1)e_3$;
\item $F_2$: $e_1\rt e_1=e_2$, $e_1\rt e_2=e_3$;
\item $F_3$: $e_1\lt e_1=e_3$, $e_1\rt e_1=e_2$, $e_1\rt e_2=e_3$;
\item $F_4(\lambda,m)$, $(\lambda,m)\in\C^2$: $e_1\lt e_1=e_1\rt e_1=e_2$, $e_2\lt e_1=\lambda e_3$, $e_2\rt e_1=(\lambda-1)e_3$, $e_1\lt e_2=e_1\rt e_2=me_3$;
\item $F_5[\lambda:m]$, $[\lambda:m]\in\Pj^1$: $e_1\lt e_1=e_1\rt e_1=e_2$, $e_2\lt e_1=e_2\rt e_1=\lambda e_3$, $e_1\lt e_2=e_1\rt e_2=me_3$;
\item $F_6[\lambda:m]$, $[\lambda:m]\in\Pj^1$: as $F_5[\lambda:m]$ but with $e_1\lt e_1=e_2+e_3$.
\end{itemize}
\emph{Algebra with two-dimensional annihilator (nilpotency index $3$):}
\begin{itemize}
\item $F_7$: $e_1\lt e_1=e_2$, $e_1\rt e_1=e_3$.
\end{itemize}
\emph{$2$-step nilpotent algebras with $S/\Ann(S)\cong\Nz$ (nilpotency index $3$):}
\begin{itemize}
\item $D(\alpha,\delta,\gamma)$, with $\alpha=(\alpha_1,\alpha_2)$, $\delta=(\delta_1,\delta_2)$, $\gamma=(\gamma_1,\gamma_2)\in\C^2$ such that $\gamma\neq0$, or $\gamma=0$ and $\alpha\neq0\neq\delta$:
\[
\begin{aligned}
&e_1\lt e_1=\alpha_1e_3,\quad e_2\lt e_2=\delta_1e_3,\quad e_1\lt e_2=-e_2\lt e_1=\gamma_1e_3,\\
&e_1\rt e_1=\alpha_2e_3,\quad e_2\rt e_2=\delta_2e_3,\quad e_1\rt e_2=-e_2\rt e_1=\gamma_2e_3;
\end{aligned}
\]
\item $T(x,y)$, with $x=(x_1,x_2)$, $y=(y_1,y_2)\in\C^2$ such that $x+y\neq0$: $e_1\lt e_2=x_1e_3$, $e_2\lt e_1=y_1e_3$, $e_1\rt e_2=x_2e_3$, $e_2\rt e_1=y_2e_3$, and moreover $e_1\lt e_1=e_3$ if $x_2+y_2\neq0$, $e_1\rt e_1=e_3$ if $x_2+y_2=0$.
\end{itemize}
Algebras with different names are not isomorphic. Within the families:
\begin{enumerate}
\item $F_1(\beta)$, $F_4(\lambda,m)$, $F_5[\lambda:m]$, $F_6[\lambda:m]$ and $\N{\alpha}{\beta}\oplus\C$ are isomorphic only for equal parameters (equal points of $\Pj^1$ for the projective ones);
\item $D(\alpha,\delta,\gamma)\cong D(\alpha',\delta',\gamma')$ if and only if there are $\lambda,\mu,\nu\in\C^*$ with $\nu^2=\lambda\mu$ such that $(\alpha',\delta',\gamma')=(\lambda\alpha,\mu\delta,\nu\gamma)$ or $(\alpha',\delta',\gamma')=(\lambda\delta,\mu\alpha,\nu\gamma)$;
\item $T(x,y)\cong T(x',y')$ if and only if $(x',y')=\nu(x,y)$ for some $\nu\in\C^*$; that is, $[x_1:x_2:y_1:y_2]\in\Pj^3$ is a complete invariant.
\end{enumerate}
\end{theorem}

Thus $D$ has four essential parameters (two when $\gamma=0$) and $T$ has three. Normal forms for $D$: if $\gamma\neq0$ one may take $\gamma=(1,\gamma_2)$ or $(0,1)$, and then $(\alpha,\delta)$ is determined up to $(t\alpha,t^{-1}\delta)$ and the swap $(\alpha,\delta)\mapsto(\delta,\alpha)$; if $\gamma=0$, the algebra depends only on the unordered pair of points $[\alpha_1:\alpha_2]$, $[\delta_1:\delta_2]\in\Pj^1$.

\section{Proof of the main theorem}

Let $S$ be a three-dimensional nilpotent LSDA. By Lemma~\ref{lem:ann}, $\dim\Ann(S)\ge1$, and $\dim\Ann(S)=3$ means $S=\mathcal N_3$. If $S$ is split, then $S=B\oplus\C$ with $B$ a two-dimensional nilpotent LSDA, so $S$ is $\mathcal N_3$ or $\N{\alpha}{\beta}\oplus\C$ by Proposition~\ref{prop:dim2}. In the latter case every $x\notin\Ann(S)$ satisfies $x\lt x=c\alpha e_2$ and $x\rt x=c\beta e_2$ with $c\neq0$, so $[\alpha:\beta]$ is an isomorphism invariant. From now on $S$ is non-split, and by Corollary~\ref{cor:proc} it corresponds to a pair $(A,W)$ with $\dim A=3-\dim\Ann(S)$.

\subsection{Two-dimensional annihilator}
If $\dim\Ann(S)=2$, then $A=\mathcal N_1$ and $W$ is a two-dimensional subspace of $H^2(\mathcal N_1,\C)=\C^2$, so $W=\langle[\Dl{11}],[\Dr{11}]\rangle$ is the whole space. The condition $\Ann(\mathcal N_1)\cap\Ann(W)=0$ holds. Taking $\theta=\Dl{11}e_2+\Dr{11}e_3$ gives $F_7$.

\subsection{One-dimensional annihilator, $A=\N{\alpha}{\beta}$}
Now $\dim\Ann(S)=1$, $W=\langle[\theta]\rangle$ and $V=\C e_3$. By Lemma~\ref{lem:iso} the relevant group is $\Aut(A)\times\C^*$ acting as in Lemma~\ref{lem:action}, on classes modulo $B^2$. Since $\Ann(A)=\langle e_2\rangle$ and $\lambda_{22}=\rho_{22}=0$, the condition $\Ann(A)\cap\Ann(\theta)=0$ means $(\lambda_{12},\lambda_{21},\rho_{12},\rho_{21})\neq0$. In each case $A_\theta$ has the products of $A$ plus the $e_3$-components given by $\theta$.

\emph{Case I: $A=\N{1}{\beta}$, $\beta\neq1$.} By Proposition~\ref{prop:H2}, $\theta=\lambda_{11}\Dl{11}+\rho_{11}\Dr{11}+s(\Dl{21}+(\beta-1)\Dr{12})$ with $s\neq0$. By Lemma~\ref{lem:action}, $s\mapsto kp^3s$ and $(\lambda_{11},\rho_{11})\mapsto k\bigl(p^2(\lambda_{11},\rho_{11})+pqs(1,\beta-1)\bigr)$. Modulo $B^2=\langle(1,\beta)\rangle$ the pair $(\lambda_{11},\rho_{11})$ is recorded by the number $\beta\lambda_{11}-\rho_{11}$, and the shift $(1,\beta-1)$ changes it by $kpqs\neq0$ per unit of $q$. Choosing $q$ makes this class zero, and $k=(p^3s)^{-1}$ gives $s=1$. Hence there is a single orbit, represented by $\theta=\Dl{21}+(\beta-1)\Dr{12}$, which gives $F_1(\beta)$.

\emph{Case II: $A=\N{0}{1}$.} Here $\theta=c\Dl{11}+s\Dr{12}$ modulo $B^2=\langle\Dr{11}\rangle$, with $s\neq0$. Since $\lambda_{12}=\lambda_{21}=0$, Lemma~\ref{lem:action} gives $(c,s)\mapsto(kp^2c,kp^3s)$. The orbits are $c=0$, giving $F_2$, and $c\neq0$ (take $p=c/s$ and $k=(p^2c)^{-1}$, so $c=s=1$), giving $F_3$. Whether $c=0$ is invariant, so $F_2\not\cong F_3$.

\emph{Case III: $A=\N{1}{1}$.} Take the representative with $\rho_{11}=0$ and write $\theta=d\Dl{11}+\lambda\Dl{21}+\rho\Dr{21}+m(\Dl{12}+\Dr{12})$, with $w=(\lambda,\rho,m)\neq0$. The class of $d=\lambda_{11}-\rho_{11}$ modulo $B^2$ is well defined, and Lemma~\ref{lem:action} gives
\[
w\mapsto kp^3w,\qquad d\mapsto k\bigl(p^2d+pq(\lambda-\rho)\bigr).
\]
If $\lambda\neq\rho$, choose $q$ with $d=0$ and $k$ with $kp^3(\lambda-\rho)=1$, so $\lambda-\rho=1$. The transformations preserving this normal form satisfy $kp^3=1$ and $q=0$, so they fix $(\lambda,m)$. This gives $F_4(\lambda,m)$, with pairwise non-isomorphic members. If $\lambda=\rho$, then $d\mapsto kp^2d$ and $w=(\lambda,\lambda,m)\mapsto kp^3w$, with $(\lambda,m)\neq0$. For $d=0$ only the point $[\lambda:m]\in\Pj^1$ survives, giving $F_5[\lambda:m]$. For $d\neq0$ put $k=(p^2d)^{-1}$, so $d=1$ and $w\mapsto(p/d)w$; again only $[\lambda:m]$ survives, giving $F_6[\lambda:m]$. The quantities $\lambda-\rho$ (zero or not) and $d$ (zero or not when $\lambda=\rho$) are invariant, so $F_4$, $F_5$ and $F_6$ are mutually non-isomorphic.

\subsection{One-dimensional annihilator, $A=\Nz$ (the $2$-step case)}
Now every pair is a cocycle and $B^2=0$, so $\theta$ is a pair of $2\times2$ matrices $B_1=(\lambda_{ij})$, $B_2=(\rho_{ij})$. Since $\Aut(\Nz)=\GL_2$, Lemma~\ref{lem:iso} says that $(B_1,B_2)$ and $(kP^{\mathsf T}B_1P,\,kP^{\mathsf T}B_2P)$ give isomorphic algebras for $P\in\GL_2$, $k\in\C^*$, and only these. The algebra is non-split if and only if there is no $x\neq0$ with $B_ix=B_i^{\mathsf T}x=0$ for $i=1,2$.

\emph{Step 1 (symmetric and skew parts).} Write $B_i=S_i+\gamma_iJ$ with $S_i$ symmetric and $J=\left(\begin{smallmatrix}0&1\\-1&0\end{smallmatrix}\right)$. Since $P^{\mathsf T}JP=\det(P)J$, the group acts by $S_i\mapsto kP^{\mathsf T}S_iP$ and $\gamma\mapsto k\det(P)\gamma$. Because $B_ix=B_i^{\mathsf T}x=0$ iff $S_ix=0$ and $\gamma_iJx=0$, the algebra is non-split iff $\gamma\neq0$, or $\gamma=0$ and $S_1,S_2$ have no common kernel vector.

\emph{Step 2 (orthogonal geometry).} Let $Q$ be the space of symmetric $2\times2$ matrices with the quadratic form $q(S)=\det S$ and its polar form $\langle\ ,\ \rangle$. The map $\rho(P)S=P^{\mathsf T}SP/\det P$ preserves $q$, its kernel is the group of scalar matrices, and its image is a connected three-dimensional subgroup of $\Ot(Q,q)$, hence equal to $\SO(Q,q)$. With $c=k\det P$ we have $kP^{\mathsf T}SP=c\,\rho(P)S$, so the action on $(S_1,S_2,\gamma)$ is $(c\rho(P)S_1,c\rho(P)S_2,c\gamma)$ with $c\in\C^*$ and $\rho(P)\in\SO(Q,q)$ independent. The non-zero isotropic vectors of $Q$ are the symmetric matrices of rank one.

\emph{Step 3 (two types).} Let $U=\langle S_1,S_2\rangle$. Either $U$ lies in a non-degenerate plane of $Q$, or $U$ is a degenerate plane; in the second case $U=\ell^\perp$ for the isotropic line $\ell=\operatorname{rad}U$, because $Q$ contains no totally isotropic plane. Every vector, and every non-degenerate plane, lies in a non-degenerate plane, and by Witt's theorem \cite{Lam} all non-degenerate planes are conjugate under $\SO(Q,q)$ to the plane $\Delta_0$ of diagonal matrices. So we get:
\begin{itemize}
\item \emph{type $D$}: after a congruence $S_1$ and $S_2$ are both diagonal, i.e.\ they are simultaneously diagonalisable by congruence;
\item \emph{type $T$}: $S_1,S_2$ are linearly independent and the binary form $\det(sS_1+tS_2)$ has a double root.
\end{itemize}
After moving $\ell$ to $\langle E_{11}\rangle$, type $T$ has both $S_i$ in $\ell^\perp=\bigl\{\left(\begin{smallmatrix}a&b\\b&0\end{smallmatrix}\right)\bigr\}$. The type is invariant under the group, so no algebra of type $D$ is isomorphic to one of type $T$.

\emph{Step 4 (type $D$).} Write $B_i=\left(\begin{smallmatrix}\alpha_i&\gamma_i\\-\gamma_i&\delta_i\end{smallmatrix}\right)$.\linebreak[3] For $P=\diag(p,r)$ the action is \[(\alpha,\delta,\gamma)\mapsto(kp^2\alpha,kr^2\delta,kpr\gamma),\] which realises exactly the triples $(\lambda,\mu,\nu)$ with $\nu^2=\lambda\mu$; for $P=J$ and $k=1$ it swaps $\alpha$ and $\delta$ and fixes $\gamma$. Conversely, let $(P,k)$ map one algebra of type $D$ in diagonal form to another, with $c=k\det P$. Then $\rho(P)$ restricts to an isometry from $\langle S_1,S_2\rangle$ onto $\langle c^{-1}S_1',c^{-1}S_2'\rangle$, both inside $\Delta_0$. By Witt's theorem in $\Delta_0$ it extends to $h\in\Ot(\Delta_0)$, and $h\oplus\det(h)$ on $\Delta_0^\perp$ equals $\rho(P_0)$ with $P_0$ diagonal or antidiagonal, since those are exactly the $P$ for which $P^{\mathsf T}\left(\begin{smallmatrix}0&1\\1&0\end{smallmatrix}\right)P$ is proportional to $\left(\begin{smallmatrix}0&1\\1&0\end{smallmatrix}\right)$. Hence the equivalence is realised by the maps above, which proves the isomorphism criterion for $D$. If $\gamma=0$, a common kernel vector exists iff $\alpha=0$ (then $e_1$) or $\delta=0$ (then $e_2$).

\emph{Step 5 (type $T$).} Since $\ell=\operatorname{rad}U$ is determined by the algebra, equivalences between normal forms fix $\langle E_{11}\rangle$, so $P^{\mathsf T}E_{11}P$ is a multiple of $E_{11}$ and $P=\left(\begin{smallmatrix}p&0\\s&r\end{smallmatrix}\right)$ is lower triangular. Writing $B_i=\left(\begin{smallmatrix}a_i&x_i\\y_i&0\end{smallmatrix}\right)$, a direct computation gives
\[
(x,y)\mapsto kpr\,(x,y),\qquad a\mapsto k\bigl(p^2a+ps\,(x+y)\bigr).
\]
Type $T$ means that $a=(a_1,a_2)$ and $x+y$ are linearly independent. Keeping $kpr$ fixed, $a$ can be moved to any vector independent of $x+y$: we take $a=(1,0)$ if $x_2+y_2\neq0$ and $a=(0,1)$ otherwise. The maps preserving this $a$ are those with $kp^2=1$ and $s=0$, and they act by $(x,y)\mapsto\nu(x,y)$ with $\nu=r/p$ arbitrary. Hence $[x_1:x_2:y_1:y_2]\in\Pj^3$ is a complete invariant. Type $T$ is always non-split: if $S_1v=S_2v=0$ with $S_i=\left(\begin{smallmatrix}a_i&b_i\\b_i&0\end{smallmatrix}\right)$ and $b\neq0$, then $b_iv_1=0$ gives $v_1=0$ and then $b_iv_2=0$ gives $v_2=0$.

This completes the proof of Theorem~\ref{thm:main}, except for the non-isomorphism of algebras with different names, which follows from Proposition~\ref{prop:noniso}. \qed

\section{Invariants, non-isomorphism and special subclasses}\label{sec:inv}

Table~\ref{tab:inv} lists invariants of all algebras of Theorem~\ref{thm:main}. Here $\Der(S)$ is the Lie algebra of derivations, i.e.\ linear maps $D$ with $D(x\circ y)=Dx\circ y+x\circ Dy$ for $\circ\in\{\lt,\rt\}$; its dimension equals $\dim\Aut(S)$. For the families, $\dim\Der$ is given at generic parameters.

\begin{table}[ht]
\caption{Invariants of the three-dimensional nilpotent LSDA.}\label{tab:inv}
\footnotesize\setlength{\tabcolsep}{4pt}
\begin{tabular}{@{}lccccccc@{}}
\toprule
Algebra & $\dim\Ann$ & $S/\Ann(S)$ & nil.\ index & $\dim S^2$ & $\dim\Der$ & $\lt=\rt$ & diassociative\\
\midrule
$F_1(\beta)$ & 1 & $\N{1}{\beta}$ & 4 & 2 & 2 & no & no\\
$F_2$ & 1 & $\N{0}{1}$ & 4 & 2 & 3 & no & no\\
$F_3$ & 1 & $\N{0}{1}$ & 4 & 2 & 2 & no & no\\
$F_4(\lambda,m)$ & 1 & $\N{1}{1}$ & 4 & 2 & 2 & no & no\\
$F_5[\lambda:m]$ & 1 & $\N{1}{1}$ & 4 & 2 & 3 & yes & iff $\lambda=m$\\
$F_6[\lambda:m]$ & 1 & $\N{1}{1}$ & 4 & 2 & 2 & no & iff $\lambda=m$\\
$F_7$ & 2 & $\mathcal N_1$ & 3 & 2 & 3 & no & yes\\
$D(\alpha,\delta,\gamma)$ & 1 & $\Nz$ & 3 & 1 & 3 & iff $\alpha_1{=}\alpha_2,\delta_1{=}\delta_2,\gamma_1{=}\gamma_2$ & yes\\
$T(x,y)$ & 1 & $\Nz$ & 3 & 1 & 3 & no & yes\\
$\N{\alpha}{\beta}\oplus\C$ & 2 & -- & 3 & 1 & 5 & iff $\alpha=\beta$ & yes\\
$\mathcal N_3$ & 3 & -- & 2 & 0 & 9 & yes & yes\\
\bottomrule
\end{tabular}
\end{table}

\begin{proposition}\label{prop:noniso}
The algebras of Theorem~\ref{thm:main} with different names are pairwise non-isomorphic.
\end{proposition}

\begin{proof}
$\dim\Ann(S)$ and the isomorphism class of $S/\Ann(S)$ are invariants. By Proposition~\ref{prop:dim2} they separate $F_1(\beta)$ (for each $\beta$), the group $\{F_2,F_3\}$, the group $\{F_4,F_5,F_6\}$, the group $\{D,T\}$, $F_7$, the split algebras and $\mathcal N_3$. Inside each group the algebras correspond to different orbits by Lemma~\ref{lem:iso} and Section~5: $c=0$ or not for $F_2,F_3$; $\lambda-\rho=0$ or not and $d=0$ or not for $F_4,F_5,F_6$; simultaneous diagonalisability for $D$ and $T$. The column $\dim\Der$ gives an independent check for $F_2\not\cong F_3$ and for $F_5\not\cong F_4,F_6$.
\end{proof}

\begin{corollary}[Left-symmetric algebras]\label{cor:lsa}
Up to isomorphism, the three-dimensional nilpotent complex left-symmetric algebras are $\mathcal N_3$, $\N{1}{1}\oplus\C$ ($e_1e_1=e_2$), $F_5[\lambda:m]$ ($e_1e_1=e_2$, $e_2e_1=\lambda e_3$, $e_1e_2=me_3$) and the $2$-step algebras $e_ie_j=B_{ij}e_3$ with
\[
B=I,\qquad B=J,\qquad B=J+E_{11},\qquad B=J+\diag(1,\mu)\ (\mu\neq0).
\]
\end{corollary}

\begin{proof}
These are the algebras of Theorem~\ref{thm:main} with $\lt=\rt$. For $D$ this means one form $B=\left(\begin{smallmatrix}a&g\\-g&d\end{smallmatrix}\right)$. If $g=0$ then $a,d\neq0$ and $B\sim I$. If $g\neq0$ we normalise $g=1$; then $(a,d)\sim(\lambda a,\lambda^{-1}d)$ and the swap, so $ad$ is the only invariant, with the two orbits $(0,0)$ and $(1,0)$ when $ad=0$. In $T$ the condition $a\parallel x+y$ would be needed, which is excluded.
\end{proof}

\begin{corollary}[Diassociative algebras]\label{cor:dias}
The three-dimensional nilpotent complex LSDA that are diassociative are $F_5[1:1]$, $F_6[1:1]$, $F_7$, all $D(\alpha,\delta,\gamma)$, all $T(x,y)$, $\N{\alpha}{\beta}\oplus\C$ and $\mathcal N_3$.
\end{corollary}

\begin{proof}
Algebras with $S^3=0$ satisfy every identity of degree three, in particular the five diassociative axioms. For $F_1$--$F_6$, the axioms on basis triples hold only for $F_5[\lambda:m]$ and $F_6[\lambda:m]$ with $\lambda=m$; for instance $(e_1\lt e_1)\lt e_1=\lambda e_3$ and $e_1\lt(e_1\lt e_1)=me_3$ in $F_5$ and $F_6$.
\end{proof}

\begin{remark}
For a generic algebra of type $D$ the stabiliser of $(B_1,B_2)$ in $\GL_2\times\C^*$ is one-dimensional, and the maps $e_i\mapsto e_i+c_ie_3$ ($i=1,2$) give two more dimensions; this agrees with $\dim\Der=3$ in Table~\ref{tab:inv}.
\end{remark}

\section{Concluding remarks}

The three-dimensional nilpotent LSDA form nine non-split classes, with a four-parameter family already in the $2$-step case. Natural next steps are the geometric classification (degenerations and irreducible components of the variety of three-dimensional nilpotent LSDA), the classification of four-dimensional nilpotent LSDA by the same method, and the classification of all three-dimensional LSDA, extending \cite[Theorem~2.6]{R16}.

\end{document}